\documentclass[preprint,11pt]{elsarticle}

\usepackage{amssymb}

\usepackage{hyperref}
\usepackage{eqnarray}
\usepackage{amsmath,amsfonts,amssymb}
\usepackage{cleveref}
\usepackage{multirow}
\usepackage{threeparttable}
\usepackage{graphicx}
\usepackage{xcolor}
\usepackage{subcaption} 
\usepackage{appendix}

\biboptions{sort&compress}

\usepackage{geometry}
\usepackage[T1]{fontenc}
\usepackage{amsthm}
\usepackage{mathrsfs}
\usepackage{mathtools}
\usepackage{upgreek}
\usepackage{dsfont}
\usepackage{enumerate}
\usepackage{enumitem}
\usepackage{xpatch}
\usepackage{calc}
\usepackage{xspace}
\usepackage[stable]{footmisc}
\usepackage{algorithm}
\usepackage{algorithmic}
\usepackage{float}
\usepackage{caption}
\usepackage{graphicx}
\usepackage{wrapfig}
\usepackage{stackengine,xcolor}
\usepackage{booktabs} 
\usepackage{pifont}

\usepackage{caption}     
\usepackage{makecell}

\crefname{section}{Appendix}{Appendices}
\Crefname{section}{Appendix}{Appendices}

\newcommand{\bv}{\boldsymbol{v}}%
\makeatletter
\newcommand{\oset}[3][0ex]{%
  \mathrel{\mathop{#3}\limits^{%
    \vbox to#1{\kern-2\ex@
    \hbox{$\scriptstyle#2$}\vss}}}}
\makeatother

\begin{document}

\begin{frontmatter}







\title{Eig-PIELM: A Mesh-Free Approach for Efficient Eigen-Analysis with Physics-Informed Extreme Learning Machines}






\author[inst1]{Rishi Mishra}
\author[inst2]{Smriti}
\author[inst1]{Ganapathy Krishnamurthi}
\author[inst3]{Balaji Srinivasan\corref{cor1}}
\author[inst2]{Sundararajan Natarajan\fnref{label2}}

\affiliation[inst1]{organization={Department of Engineering Design},
            addressline={Indian Institute of Technology Madras}, 
            city={Chennai},
            postcode={600036}, 
            state={Tamil Nadu},
            country={India}}

\affiliation[inst2]{organization={Department of Mechanical Engineering},
            addressline={Indian Institute of Technology Madras}, 
            city={Chennai},
            postcode={600036}, 
            state={Tamil Nadu},
            country={India}}

\affiliation[inst3]{organization={Wadhwani School of Data Science \& AI},
            addressline={Indian Institute of Technology Madras}, 
            city={Chennai},
            postcode={600036}, 
            state={Tamil Nadu},
            country={India}}

\cortext[cor1]{Corresponding author}
\fntext[label2]{Pandurangan Faculty Fellow}

\begin{abstract}
In this work, a novel Eig-PIELM framework is proposed that extends physics-informed extreme learning machine for an efficient and accurate solution of linear eigenvalue problems. The method reformulates the governing differential equations into a compact algebraic system solvable in a single step. Boundary conditions are enforced exactly via an algebraic projection onto the boundary-admissible subspace, eliminating the computational overhead of penalty parameters, and backpropagation while preserving the computational advantages of extreme learning machines. The proposed framework is mesh-free and yields both eigenvalues and mode shapes simultaneously in one linear solve. The robustness and accuracy of the proposed framework is demonstrated through a range of benchmark problems. We believe that the mesh-free nature, solution structure and accuracy of Eig-PIELM makes it particularly valuable for parametric studies in mechanical, acoustic, and electromechanical systems where rapid frequency spectrum analysis is critical.

\end{abstract}



\end{frontmatter}


\section{Introduction}
\label{sec:intro}
\noindent The numerical solution of eigenvalue problems represents a fundamental challenge across computational science and engineering. These problems seek non-trivial pairs $(\lambda,\boldsymbol v)$ such that a linear operator acting on a function or vector reproduces that object up to a scalar factor. In finite-dimensional linear algebra this reads:
\begin{equation}
    \mathbf{A}\boldsymbol v=\lambda  \boldsymbol v, \qquad \text{with} \qquad \boldsymbol v \neq \boldsymbol{0},
\end{equation}
with $\mathbf{A}\in\mathbb R^{n\times n}$ or $\mathbb C^{n\times n}$ being the system matrix and $\bv$ the eigenvector. In the continuous setting, one encounters the Sturm–Liouville or more general partial-differential form:
\begin{equation}
    \mathcal{L} \boldsymbol v=\lambda \boldsymbol v, \quad \boldsymbol v ~~ \text{in} ~~\Omega,\quad\text{with}\qquad \mathcal{B} \boldsymbol v=0 ~~\text{on}~~\partial \Omega.
\end{equation}
with $\mathcal{L}$ being the differential operator, $\mathcal{B}$ the boundary operator, $\Omega$ the spatial domain and $\partial\Omega$ its boundary. These problems are inherently non-linear since both the eigenvalues, $\lambda$ and the eigenfunctions, $\boldsymbol v$ are unknowns. Such eigen-analyses reveal natural modes in quantum mechanics \cite{schrodinger1940method,turbiner1984eigenvalue}, structural vibration and buckling \cite{bathe1973solution,shinozuka1972random}, fluid-mechanical stability \cite{cliffe1994eigenvalues,mahajan1991eigenvalue}, and electromagnetic or photonic resonance \cite{valovik2018nonlinear,guo2009eigen}.\\
\newline
Classical eigenvalue computation begins with solving the characteristic equation $\det(\mathbf{A} - \lambda \mathbf{I}) = 0$ for small matrices and exploiting spectral decomposition for symmetric/Hermitian cases \cite{strang1993linear,horn2013matrix}. When dealing with continuum problems, these are typically discretized into an algebraic form through finite differences \cite{choi1986finite,baxley1972eigenvalues}, finite element assembly of sparse stiffness and mass matrices \cite{bertrand2023reduced,hughes2014finite,chen2011finite}, or spectral methods \cite{feit1982solution, gheorghiu2014spectral, canuto2007spectral}. The resulting discrete eigenvalue problems are then solved using the QR algorithm with orthogonal transformations \cite{Lanczos1950,Francis1961}, iterative methods like power iteration \cite{parlett1998symmetric,trefethen1997numerical}, or Krylov subspace approaches \cite{Saad1980} implemented in libraries such as ARPACK \cite{Lehoucq1998}, SLEPc, and ELPA \cite{saad2011numerical}.\\
\newline
Recent advances have introduced neural network methods for eigenvalue problems. Lagaris and Fotiadis (1997) first introduced multilayer perceptrons with ansatz-based boundary enforcement for solving the Schrödinger equation \cite{lagaris1997artificial}. This was subsequently extended by Yu and E’s Deep Ritz method \cite{yu2018deep}. Raissi et al. \cite{raissi2019physics} formalized PINNs, that penalize PDE and boundary residuals via automatic differentiation. For eigenproblems, PINNs optimize eigenfunctions and eigenvalues jointly (or via Rayleigh quotients) with normalization/orthogonality constraints \cite{jin2020unsupervised,jin2022physics,ben2023deep}. Despite successes across physics domains \cite{yang2023data,yu2024solving,zhang2024orthogonal,pallikarakis2024application,singhal2024physics,mattheakis2022first,chen2024pinn,yoo2025physics}, PINN training can suffer from nonconvexity and ill-conditioning, sensitivity to initialization/hyperparameters, high cost in low dimensions, and spectral bias against high-frequency modes \cite{krishnapriyan2021characterizing,wang2022understanding,grossmann2023can,rahaman2019spectral,tancik2020fourier,wang2021eigenvector}.\\
\newline
Physics-Informed Extreme Learning Machines (PIELMs) address training cost by freezing a single hidden layer and solving for output weights in one linear least-squares step (via Moore–Penrose or regularized normal equations), embedding differential and boundary operators directly \cite{dwivedi2020physics}. This has been applied to a range of PDEs \cite{pan2024applying,liu2023bayesian,wang2025physics,li2023augmented,yan2022framework,joshi2024physics,ren2025physics}. However, a vanilla PIELM faces a core obstacle for eigenproblems: the residual depends on both $\lambda$ and the coefficients, breaking the single linear solve; existing work either assumes $\lambda$ is known \cite{huang2025physics} or resorts to iteration, losing PIELM’s efficiency.\\
\newline
In this paper, we propose \textbf{Eig-PIELM}, which retains PIELM’s direct-solve character while handling eigenproblems. The key idea is an algebraic projection that restricts coefficients to a boundary-admissible subspace, satisfying $\mathcal{B}\boldsymbol v=0$ exactly at boundary collocation points and leaving only interior residuals. Stationarity with respect to the output coefficients and the eigenvalue yields a coupled system that collapses to a \emph{single} symmetric generalized eigenproblem solvable by standard routines. The method is mesh-free and non-iterative, computes spectra and modes in one solve, and avoids penalty weights and backpropagation. Numerical tests in structural vibrations, acoustics, and elastodynamics show high accuracy at very low computational cost. The main contributions of this work include:
\begin{itemize}
       \item[\ding{113}] A novel Eig-PIELM framework that converts differential eigenproblems into a single symmetric generalized eigensystem.
        \item[\ding{113}] Exact boundary enforcement via algebraic projection, eliminating penalty weights and tuning.
        \item[\ding{113}] One-shot, mesh-free computation of multiple eigenpairs without iterative optimization.
        \item[\ding{113}] Extensive benchmarks (vibrations, acoustics) demonstrating accuracy and speed.
\end{itemize}
The paper proceeds with the mathematical setup in Section~\ref{sec:goveq}; Section~\ref{sec:pielmtheory} details PIELM and the proposed Eig-PIELM framework. Section~\ref{sec:numex} presents results followed by conclusions in the last section.

\section{Governing differential equations for vibrating systems}\label{sec:goveq}
\noindent In this section, the governing equations for vibrating systems considered in this study are summarized, these include, axial vibration of a bar, transverse vibration of a beam, the Helmholtz equation for wave propagation, and free vibration of a two-dimensional elastic body. In all the examples considered in this study, the domain is assumed to be homogeneous and follows linear elastic theory where applicable.


\subsection{Euler-Bernoulli Beam}\label{eig_gov_eqn}
\noindent Within the framework of the Euler-Bernoulli beam theory, the governing equation is given by
\begin{equation}
\left\{
\begin{aligned}
   & \dfrac{\partial^2}{\partial x^2} \left( EI \dfrac{\partial^2 W}{\partial x^2} \right)+ \rho A\dfrac{\partial^2 W}{\partial t^2}= 0 \quad {\rm in} \quad \Omega, \\
   &W=g_{D_{w}} \quad {\rm and} \quad
    \dfrac{\partial W}{\partial x}=g_{D_{\theta}} \quad  ~~~~{\rm on} \quad {\partial\Omega}_{D},\\
    &\dfrac{\partial^2 W}{\partial x^2} = g_{M} \quad  {\rm and}  \quad
    \dfrac{\partial^3 W}{\partial x^3} =g_{v} \quad ~~ {\rm on} \quad {\partial\Omega}_{N},
    \label{eqn:ebt_govern}
\end{aligned} \right.
\end{equation}
where $f$ is the external excitation, $E$ is the Young's modulus, $I$ the moment of inertia, $EI$ is the flexural rigidity of the beam, $\rho$ is the material density, and $A$ is the cross-section area and 
${\partial\Omega}_{D}$ and ${\partial\Omega}_{N}$ are the regions where Dirichlet and Neumann boundary conditions are specified, respectively. Assuming harmonic solution, i.e., $W(x,t)=w(x)e^{-i\omega t}$ in \Cref{eqn:ebt_govern} and in the absence of external excitation, the governing equation for free vibration of axially functionally graded Euler-Bernoulli beam is given by
\begin{equation}
    \dfrac{\partial^2}{\partial x^2} \left( EI \dfrac{\partial^2 w}{\partial x^2} \right)- \rho A \omega^2 w= 0 \quad {\rm in} \quad \Omega \, \qquad\, {\rm with~Boundary~conditions,}
    \label{eqn:ebt_beam_freq}
\end{equation}
where $\omega$ is the eigenfrequency.

\subsection{Acoustic wave equation}
\noindent The wave propagation in a linear acoustic medium is governed by the following scalar equation:
\begin{equation}
    \dfrac{1}{c^2}\,\dfrac{\partial^2 p}{\partial t^2} = \nabla^2 p, \qquad {\rm in} \quad \Omega,
    \label{eqn:helmholtz}
\end{equation}
where $p$ and $c$ represents the acoustic pressure and the speed of sound, respectively and the gradient operator $\nabla$ is formulated in Cartesian coordinates as $\nabla^\top = \left[ \frac{\partial}{\partial x} \,\, \frac{\partial}{\partial y} \right]^\top$. Upon assuming a separable harmonic form $p(x,y,t) = P(x,y) e^{-i\omega t}$, we get the following governing equation:
\begin{equation}
\left\{
\begin{aligned}
&\nabla^2 P-k^2 P = 0, \quad (x,y) \in \Omega, \\
&p = \hat{p} \quad \textup{on} \quad {\partial\Omega}_D, \\
&\nabla p \cdot \mathbf{n} = 0 \quad \textup{on} \quad {\partial\Omega}_N,
\end{aligned}
\right.
\label{eqn:helm_freq}
\end{equation}
where $k=\omega/c$ is the wave number and $\omega$ is the angular frequency in rad/s.

\section{PIELM}
\label{sec:pielmtheory}
\subsection{Brief Review of Physics-Informed Extreme Learning Machine (PIELM)}
\noindent PIELM uses a single hidden layer with fixed (frozen) hidden features, as in ELM, and determines only the output coefficients by minimizing the physics residual in a least-squares sense, as in PINNs. Training therefore reduces to a single linear solve without backpropagation or iterative optimization. Let $\Omega\subset\mathbb{R}^d$ be the domain with boundary $\partial\Omega$ and $\bar\Omega=\Omega\cup\partial\Omega$. Given a linear differential operator $\mathcal{L}$ of order $m$ and a boundary operator $\mathcal{B}$, we seek $u\in C^{m}(\bar\Omega)$ satisfying
\begin{equation}
    \mathcal{L}u(\mathbf{x})=f(\mathbf{x})\quad(\mathbf{x}\in\Omega),\qquad
  \mathcal{B}u(\mathbf{x}^b)=g(\mathbf{x}^b)\quad(\mathbf{x}^b\in\partial\Omega).
\end{equation}
Next, we choose an activation $\sigma\in C^m(\mathbb{R})$ and fix $N_\phi$ features (random or polynomial). The PIELM ansatz is
\begin{equation}
     \widehat u(\mathbf{x};\boldsymbol\beta)
  =\sum_{k=1}^{N_\phi}\beta_k\,\sigma(\mathbf{w}_k^{\!\top}\mathbf{x}+b_k)
  =\boldsymbol\phi(\mathbf{x})^{\!\top}\boldsymbol\beta,
\end{equation}
where $\boldsymbol\phi(\mathbf{x})=[\sigma(\mathbf{w}_1^{\!\top}\mathbf{x}+b_1),\dots,\sigma(\mathbf{w}_{N_\phi}^{\!\top}\mathbf{x}+b_{N_\phi})]^{\!\top}$ and $\boldsymbol\beta\in\mathbb{R}^{N_\phi}$ are the unknown output coefficients.\\
\newline 
On distributing interior points $\{\mathbf{x}_i\}_{i=1}^{N_x}\subset\Omega$ and boundary points $\{\mathbf{x}^b_j\}_{j=1}^{N_b}\subset\partial\Omega$, the pointwise residuals can be defined as
\begin{equation}
     r_i^{\mathrm{PDE}}(\boldsymbol\beta)=\mathcal{L}\widehat u(\mathbf{x}_i;\boldsymbol\beta)-f(\mathbf{x}_i),
  \qquad
  r_j^{\mathrm{BC}}(\boldsymbol\beta)=\mathcal{B}\widehat u(\mathbf{x}^b_j;\boldsymbol\beta)-g(\mathbf{x}^b_j).
\end{equation}
The operator actions on the fixed features are collected into
\begin{equation}
     \boldsymbol\Psi:=\begin{bmatrix}(\mathcal{L}\boldsymbol\phi(\mathbf{x}_1))^{\!\top}\\ \vdots\\ (\mathcal{L}\boldsymbol\phi(\mathbf{x}_{N_x}))^{\!\top}\end{bmatrix}\!\in\mathbb{R}^{N_x\times N_\phi},
  \qquad
  \mathbf{B}_{bc}:=\begin{bmatrix}(\mathcal{B}\boldsymbol\phi(\mathbf{x}^b_1))^{\!\top}\\ \vdots\\ (\mathcal{B}\boldsymbol\phi(\mathbf{x}^b_{N_b}))^{\!\top}\end{bmatrix}\!\in\mathbb{R}^{N_b\times N_\phi},
\end{equation}
and the targets $\mathbf{f}=[f(\mathbf{x}_1),\dots,f(\mathbf{x}_{N_x})]^{\!\top}$, $\mathbf{g}=[g(\mathbf{x}^b_1),\dots,g(\mathbf{x}^b_{N_b})]^{\!\top}$. Stacking interior and boundary rows yields
\begin{equation}
     \boldsymbol{\mathcal{X}}:=\begin{bmatrix}\boldsymbol\Psi\\ \mathbf{B}_{bc}\end{bmatrix}\in\mathbb{R}^{(N_x+N_b)\times N_\phi},
  \qquad
  \boldsymbol{\mathcal{Y}}:=\begin{bmatrix}\mathbf{f}\\ \mathbf{g}\end{bmatrix}\in\mathbb{R}^{N_x+N_b}.
\end{equation}
The physics-informed least-squares fit
\begin{equation}
    J(\boldsymbol\beta)=\tfrac12\|\boldsymbol{\mathcal{X}}\boldsymbol\beta-\boldsymbol{\mathcal{Y}}\|_2^2
\end{equation}
leads to the normal equations
\begin{equation}
   \boldsymbol{\mathcal{X}}^{\!\top}\boldsymbol{\mathcal{X}}\,\boldsymbol\beta
  =\boldsymbol{\mathcal{X}}^{\!\top}\boldsymbol{\mathcal{Y}}
  \;\;\Longleftrightarrow\;\;
  \mathbf{H}\boldsymbol\beta=\mathbf{q}, 
\end{equation}
with Gram matrix $\mathbf{H}:=\boldsymbol\Psi^{\!\top}\boldsymbol\Psi+\mathbf{B}_{bc}^{\!\top}\mathbf{B}_{bc}$ and right-hand side $\mathbf{q}:=\boldsymbol\Psi^{\!\top}\mathbf{f}+\mathbf{B}_{bc}^{\!\top}\mathbf{g}$. When $\mathbf{H}$ is nonsingular, $\boldsymbol\beta^*=\mathbf{H}^{-1}\mathbf{q}$; otherwise use a regularized or pseudoinverse solve, e.g.\ $(\mathbf{H}+\varepsilon\mathbf{I})^{-1}\mathbf{q}$ or $\mathbf{H}^{\dagger}\mathbf{q}$. The approximation $\widehat u(\mathbf{x};\boldsymbol\beta^*)$ satisfies the PDE and boundary data in the least-squares sense.

\subsection{Proposed Eig-PIELM framework for eigenvalue problems}
\label{sec:eig-pielm}
\noindent Let $\Omega\subset\mathbb{R}^d$ be a bounded domain with Lipschitz boundary $\partial\Omega$, $\bar\Omega=\Omega\cup\partial\Omega$, and let $\mathcal{L}$ be a linear differential operator of order $m$. We seek nontrivial eigenpairs $(\lambda,u)$ with $u\in C^{m}(\bar\Omega)$ such that
\begin{equation}
  \mathcal{L}\,u(\mathbf{x})=\lambda\,u(\mathbf{x})\quad(\mathbf{x}\in\Omega),\qquad
  \mathcal{B}\,u(\mathbf{x}^b)=0\quad(\mathbf{x}^b\in\partial\Omega).
  \label{eq:eig-cont}
\end{equation}
Let $\{\mathbf{x}_i\}_{i=1}^{N_x}\subset\Omega$ and $\{\mathbf{x}^b_j\}_{j=1}^{N_b}\subset\partial\Omega$ be interior and boundary collocation sets. We approximate $u$ by a single-hidden-layer PIELM expansion with frozen features (random or polynomial):
\begin{equation}
  \widehat u(\mathbf{x};\boldsymbol\beta)=\sum_{k=1}^{N_\phi}\beta_k\,\sigma(\mathbf{w}_k\!\cdot\!\mathbf{x}+b_k)
  =\boldsymbol\phi(\mathbf{x})^{\!\top}\boldsymbol\beta,
  \label{eq:ansatz}
\end{equation}
where $\sigma\in C^{m}(\mathbb{R})$, $\boldsymbol\phi(\mathbf{x})\in\mathbb{R}^{N_\phi}$ collects the basis functions, and $\boldsymbol\beta\in\mathbb{R}^{N_\phi}$ are the trainable coefficients ($N_\phi$ is the number of basis functions). Define the interior residuals $r_i^{\text{PDE}}(\boldsymbol\beta,\lambda):=\mathcal{L}\widehat u(\mathbf{x}_i)-\lambda\,\widehat u(\mathbf{x}_i)$. Using only the interior residuals, the least-squares loss expands to the quadratic form
\begin{equation}
  J(\boldsymbol\beta,\lambda)
  =\tfrac12\,\boldsymbol\beta^{\!\top}\bigl(\boldsymbol{\mathcal A}-\lambda\,\boldsymbol{\mathcal P}+\lambda^2\boldsymbol{\mathcal G}\bigr)\boldsymbol\beta,
  \label{eq:loss}
\end{equation}
with design matrices
\begin{equation}
  \boldsymbol{\mathcal A}:=\sum_{i=1}^{N_x}(\mathcal{L}\boldsymbol\phi_i)(\mathcal{L}\boldsymbol\phi_i)^{\!\top},\quad
  \boldsymbol{\mathcal S}:=\sum_{i=1}^{N_x}(\mathcal{L}\boldsymbol\phi_i)\boldsymbol\phi_i^{\!\top},\quad
  \boldsymbol{\mathcal P}:=\boldsymbol{\mathcal S}+\boldsymbol{\mathcal S}^{\!\top},\quad
  \boldsymbol{\mathcal G}:=\sum_{i=1}^{N_x}\boldsymbol\phi_i\boldsymbol\phi_i^{\!\top},
  \label{eq:design}
\end{equation}
where $\boldsymbol\phi_i:=\boldsymbol\phi(\mathbf{x}_i)$ and $\mathcal{L}\boldsymbol\phi_i:=(\mathcal{L}\phi_1(\mathbf{x}_i),\ldots,\mathcal{L}\phi_{N_\phi}(\mathbf{x}_i))^{\!\top}$. Stationarity of \eqref{eq:loss} with respect to $\boldsymbol\beta$ and $\lambda$ gives the coupled system
\begin{equation}
\label{eq:coupled}
\begin{cases}
(\boldsymbol{\mathcal A}-\lambda\,\boldsymbol{\mathcal P}+\lambda^2\boldsymbol{\mathcal G})\,\boldsymbol\beta=\mathbf{0},\\[3pt]
\tfrac12\,\boldsymbol{\mathcal P}\,\boldsymbol\beta=\lambda\,\boldsymbol{\mathcal G}\,\boldsymbol\beta,
\end{cases}
\end{equation}
where the first line is a quadratic eigenvalue problem in $\lambda$ and the second is linear in $\lambda$.\\
\newline
To obtain a compact symmetric eigenproblem that satisfies the boundary conditions exactly at the boundary collocation set, we introduce a full-column linear map $\mathsf{T}_{\!b}$ (depending on $\mathcal{B}$, the basis, and the boundary collocation) and reparameterize the coefficients as $\boldsymbol\beta=\mathsf{T}_{\!b}\mathbf{y}$. This map places the search in a boundary-admissible coordinate space. Applying this change of variables to the quadratic forms defines the reduced matrices
\[
  \boldsymbol{\mathcal G}_{\mathrm{red}}:=\mathsf{T}_{\!b}^{\!\top}\boldsymbol{\mathcal G}\,\mathsf{T}_{\!b},\qquad
  \boldsymbol{\mathcal S}_{\mathrm{red}}:=\mathsf{T}_{\!b}^{\!\top}\boldsymbol{\mathcal S}\,\mathsf{T}_{\!b},\qquad
  \boldsymbol{\mathcal P}_{\mathrm{red}}:=\mathsf{T}_{\!b}^{\!\top}\boldsymbol{\mathcal P}\,\mathsf{T}_{\!b}.
\]
Because the boundary is satisfied exactly in this coordinate space, the boundary-induced skew part disappears and $\boldsymbol{\mathcal S}_{\mathrm{red}}$ is symmetric; hence
\[
  \boldsymbol{\mathcal P}_{\mathrm{red}}=\boldsymbol{\mathcal S}_{\mathrm{red}}+\boldsymbol{\mathcal S}_{\mathrm{red}}^{\!\top}=2\,\boldsymbol{\mathcal S}_{\mathrm{red}}.
\]
Consequently, the reduced counterpart of \eqref{eq:coupled} collapses to the single symmetric generalized eigenvalue problem
\begin{equation}
\label{eq:final-gep}
\boxed{\;\boldsymbol{\mathcal S}_{\mathrm{red}}\,\mathbf{y}=\lambda\,\boldsymbol{\mathcal G}_{\mathrm{red}}\,\mathbf{y}\;}
\end{equation}
with $\boldsymbol{\mathcal G}_{\mathrm{red}}$ symmetric positive definite. Standard eigensolvers provide $(\lambda,\mathbf{y})$, and the mode reconstructs as
\[
  \widehat u(\mathbf{x})=\boldsymbol\phi(\mathbf{x})^{\!\top}\,\mathsf{T}_{\!b}\mathbf{y}.
\]

\section{Numerical Experiments}
\label{sec:numex}
\noindent To evaluate the performance of Eig-PIELM, we conduct experiments on structural dynamics and wave propagation applications. The test cases include the following free vibration problems:
\begin{enumerate}


\item Transverse vibrations of Euler-Bernoulli beams with simply supported, fixed-fixed, and fixed-free boundary configurations, involving fourth-order differential operators that challenge the method's ability to handle higher-order derivatives and mixed boundary condition.

\item Homogeneous Helmholtz equation on a rectangular domain with homogeneous Neumann boundary conditions. This is a canonical two-dimensional elliptic eigenproblem with a doubly indexed analytical spectrum.


\end{enumerate}
 For both of the cases, we assess accuracy using the absolute error $\Delta\omega_k = |\omega_k^{\text{Eig-PIELM}} - \omega_k^{\text{Exact}}|$, where $k$ is the mode number. We compute the first five eigenvalues and compare predicted mode shapes against exact solutions.\\
\newline
We employ Bernstein polynomials as basis functions within the extreme learning machine framework, chosen for their favorable numerical properties including non-negative partition of unity, endpoint interpolation capabilities, and excellent conditioning characteristics. The neural network architecture comprises $n+1$ neurons, where $n$ represents the polynomial degree, ensuring exact representation of polynomial functions up to degree $n$. Interior and boundary collocation points are sampled uniformly throughout their respective domains to ensure consistent enforcement of governing equations and boundary conditions without spatial bias.\\
\newline
All computations are performed in MATLAB R2023b using double precision arithmetic on an Intel Core i9-13980HX (2.20 GHz) processor  with 32GB RAM, without multithreading or GPU acceleration. The experimental design focuses on the first five eigenmodes for each configuration, as these lower-order modes are typically most significant for engineering applications.
\subsection{Transverse vibration of Euler–Bernoulli beam}
\label{sec:fixed_fixed_transverse}
\noindent In this section, we study the transverse vibration of beams using Eig-PIELM. The beam kinematics is governed by Euler-Bernoulli beam theory (c.f. \Cref{eqn:ebt_govern}). The beam is assumed to be isotropic and homogeneous. The material properties are taken as Young’s modulus $E = 210\,\text{GPa}$ and mass density $\rho = 7800\,\text{kg/m}^3$. The beam has length $L = 0.12\,\text{m}$, width $b = 0.003\,\text{m}$, and height $h = 0.002\,\text{m}$. The second moment of area is $I = b h^3 / 12$ and the cross-sectional area is $A = b h$. Three different boundary conditions are studied, viz., simply supported, fixed-fixed and cantilever beam. In all the cases, the results from the present framework are compared with analytical solution.

\begin{figure}
    \centering
    \includegraphics[width=0.495\linewidth,trim=30 19 21 18,clip]{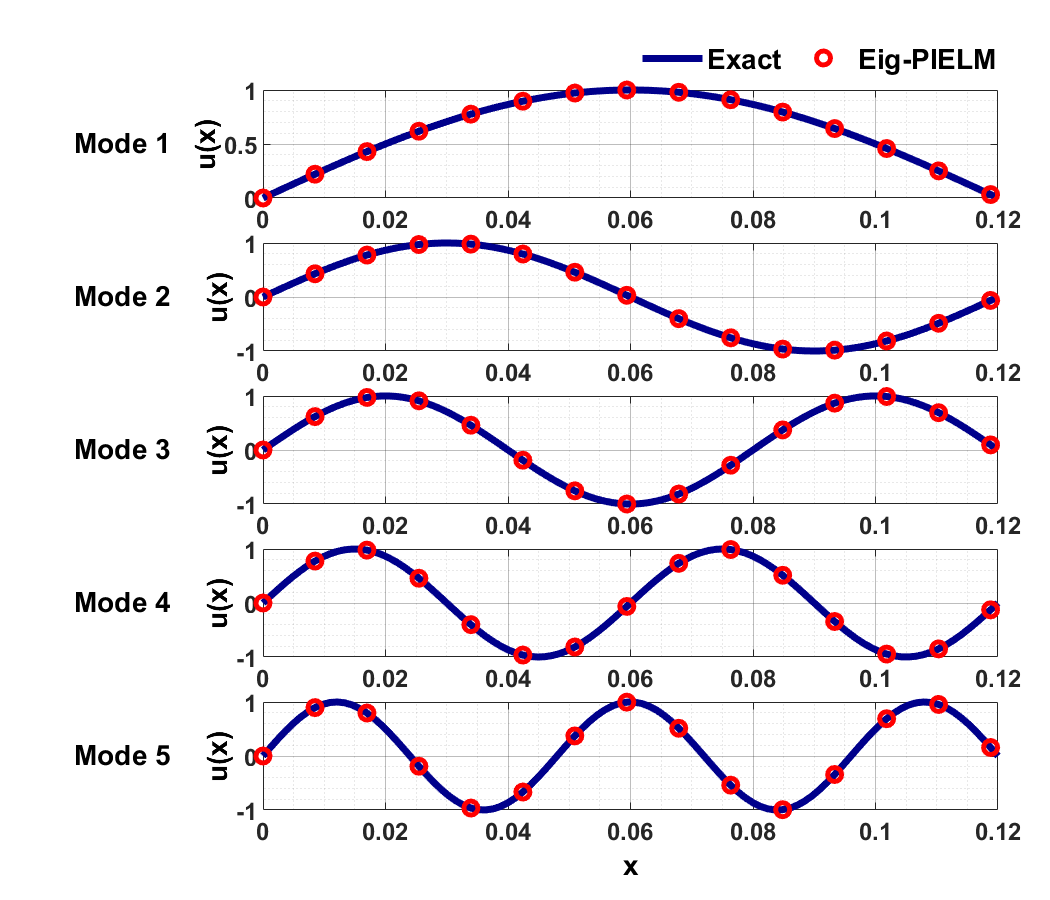}
    \caption{First five transverse vibration mode shapes of a one-dimensional simply supported Euler–Bernoulli beam from Eig-PIELM and exact analytical solutions. Mode shapes are normalized and shown over the span $x \in [0,\,L]$.}
    \label{fig:1d_bending_simply_supported}
\end{figure}

\paragraph{Simply supported} This case considers a uniform Euler–Bernoulli beam with simply supported ends. 
For a simply supported beam, the deflection vanishes at both ends, and the bending moment is zero:
\begin{equation}
\left\{ 
\begin{aligned}
&w(0) = 0, \quad \frac{d^2w}{dx^2}(0) = 0, \\
    &w(L) = 0, \quad \frac{d^2w}{dx^2}(L) = 0.
\end{aligned}
\right.
\end{equation}
The exact natural frequencies and corresponding normalized mode shapes \cite{rao2007vibration} are given by:
\begin{equation}
\begin{aligned}
    &\omega_n = \left( \frac{n\pi}{L} \right)^2
    \sqrt{\frac{E I}{\rho A}}, \\
    &\phi_m(x) = \sin\!\left( \frac{n\pi x}{L} \right)
\end{aligned}
\end{equation}
where $n=1,2,\dots.$. \Cref{fig:1d_bending_simply_supported} illustrates the first five mode shapes predicted by Eig-PIELM together with the exact analytical solutions. \Cref{tab:1d_transverse_simply_supported} compares the predicted and analytical eigenfrequencies.

\begin{table}[!ht]
\centering
\normalsize
\begin{threeparttable}
\captionsetup{labelformat=empty, font=normalsize}
\caption{\textbf{Table \thetable} \\ Comparison of predicted and analytical eigenfrequencies for the first five transverse vibration modes of a one-dimensional simply supported Euler–Bernoulli beam.}
\label{tab:1d_transverse_simply_supported}
\begin{tabular}{lc@{\hspace{2em}}c@{\hspace{2em}}c@{\hspace{2em}}c@{\hspace{2em}}c}
\toprule
Mode & $\omega_{\scriptscriptstyle\text{NP-PIELM}}$ & $\omega_{\scriptscriptstyle\text{Exact}}$ & $|\omega_{\scriptscriptstyle\text{NP-PIELM}}-\omega_{\scriptscriptstyle\text{Exact}}|$ \\
Number & $\scriptstyle(\text{rad/s})$ & $\scriptstyle(\text{rad/s})$ & $\scriptstyle(\text{rad/s})$ \\
\midrule
 1 &   2.0532$\times$10$^{3}$ &   2.0532$\times$10$^{3}$ &  1.3188$\times$10$^{-11}$ \\ 
 2 &   8.2129$\times$10$^{3}$ &   8.2129$\times$10$^{3}$ &  1.0914$\times$10$^{-11}$ \\ 
 3 &   1.8479$\times$10$^{4}$ &   1.8479$\times$10$^{4}$ &  1.0000$\times$10$^{-16}$ \\ 
 4 &   3.2852$\times$10$^{4}$ &   3.2852$\times$10$^{4}$ &  1.4552$\times$10$^{-11}$ \\ 
 5 &   5.1331$\times$10$^{4}$ &   5.1331$\times$10$^{4}$ &  8.8767$\times$10$^{-10}$ \\ 
\hline
\end{tabular}
\end{threeparttable}
\end{table}
\noindent
From Table~\ref{tab:1d_transverse_simply_supported}, it is evident that the predicted frequencies match the exact analytical values to near machine precision.  
Absolute differences are between $1.00\times 10^{-16}$ and $8.87\times 10^{-10}\,\mathrm{rad/s}$, corresponding to relative errors on the order of $10^{-20}$ to $10^{-14}$. Both low-order and higher-order modes are reproduced with extremely high accuracy. The computation, using $N_\phi = 24$ neurons and $N_x = 50{,}000$ interior points, was completed in $0.199\,\mathrm{s}$. The uniformly small errors across all modes indicate that Eig-PIELM is well-suited for accurately solving transverse vibration problems in simply supported beam configurations.

\paragraph{Fixed-Fixed} In this case, the boundary conditions are:
\begin{equation}
\left\{
\begin{aligned}
    &w(0) = 0, \quad \frac{dw}{dx}(0) = 0, \\
    &w(L) = 0, \quad \frac{dw}{dx}(L) = 0.
\end{aligned}
\right.
\end{equation}
The exact natural frequencies \cite{rao2007vibration} are given by
\begin{equation}
    \omega_n = \left( \frac{\beta_n}{L} \right)^2
    \sqrt{\frac{E I}{\rho A}}, 
    \quad n = 1, 2, \dots,
\end{equation}
where $\beta_n = \{4.73,\, 7.8532,\, 10.9956,\, 14.1372,\, 17.2787\}$ are the dimensionless wavenumbers for the first five fixed–fixed modes and the corresponding normalized mode shapes \cite{rao2007vibration} are given by
\begin{equation}
    \phi_n(x) = \left[\cos\!\left( \frac{\beta_n x}{L} \right) - \cosh\!\left( \frac{\beta_n x}{L} \right)\right] 
    - \left( \dfrac{\cos\beta_n - \cosh\beta_n}{\sin\beta_n - \sinh\beta_n}\right) \left[ \sin\!\left( \frac{\beta_n x}{L} \right) -\sinh\!\left( \frac{\beta_n x}{L} \right) \right],
\end{equation}

\begin{figure}[htbp]
  \centering
  \begin{subfigure}[t]{0.495\textwidth}
    \centering
    \includegraphics[width=\linewidth,trim=30 19 21 18,clip]{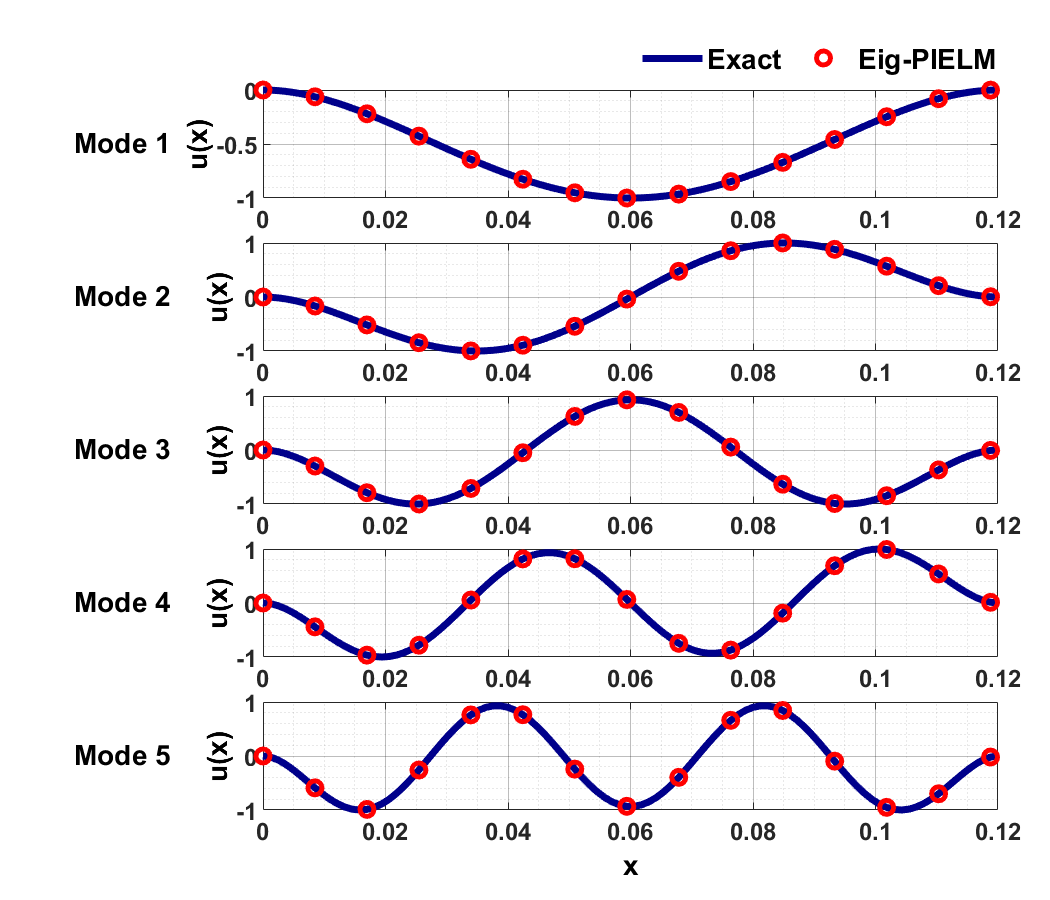}
    \caption{}
    \label{fig:1d_transverse_fixed_fixed}  
  \end{subfigure}
  \hfill
  \begin{subfigure}[t]{0.495\textwidth}
    \centering
    \includegraphics[width=\linewidth,trim=30 19 21 18,clip]{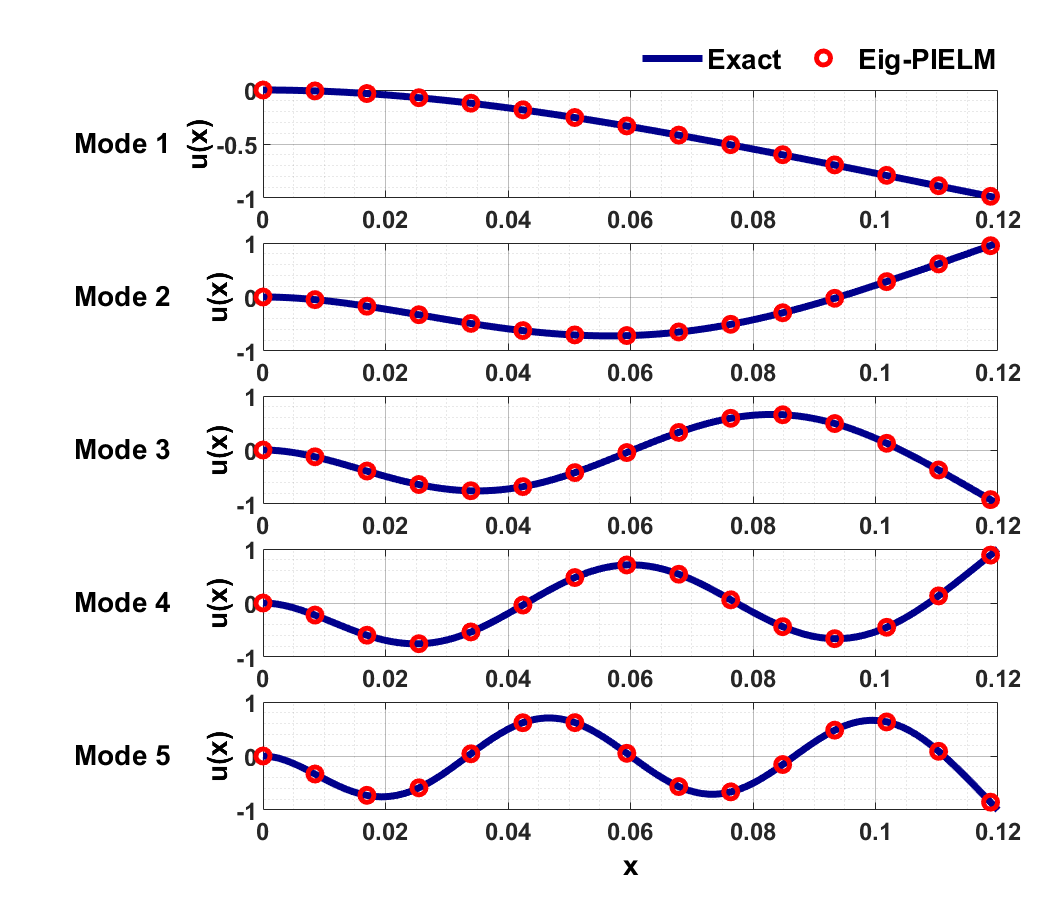}
    \caption{}
    \label{fig:1d_transverse_fixed_free}  
  \end{subfigure}
  \caption{First five transverse vibration mode shapes of a one-dimensional Euler–Bernoulli beam from Eig-PIELM and exact solutions: (a) fixed–fixed boundary conditions, and (b) fixed–free boundary conditions. The mode shapes are normalized and shown over the domain $x \in [0,\,L]$.}
\end{figure}

\noindent
\Cref{fig:1d_transverse_fixed_fixed} shows the first five mode shapes, and Table~\ref{tab:1d_transverse_fixed} compares the eigenfrequencies predicted by Eig-PIELM with the exact analytical values.
\begin{table}[!ht]
\centering
\normalsize
\begin{threeparttable}
\captionsetup{labelformat=empty, font=small}
\caption{\textbf{Table \thetable} \\ Comparison of predicted and analytical eigenfrequencies for the first five transverse vibration modes of a one-dimensional Euler–Bernoulli beam with fixed–fixed end conditions.}
\label{tab:1d_transverse_fixed}
\begin{tabular}{lc@{\hspace{2em}}c@{\hspace{2em}}c@{\hspace{2em}}c@{\hspace{2em}}c}
\toprule
Mode & $\omega_{\scriptscriptstyle\text{NP-PIELM}}$ & $\omega_{\scriptscriptstyle\text{Exact}}$ & $|\omega_{\scriptscriptstyle\text{NP-PIELM}}-\omega_{\scriptscriptstyle\text{Exact}}|$ \\
Number & $\scriptstyle(\text{rad/s})$ & $\scriptstyle(\text{rad/s})$ & $\scriptstyle(\text{rad/s})$ \\
\midrule
 1 &   4.6545$\times$10$^{3}$ &   4.6545$\times$10$^{3}$ &  1.8190$\times$10$^{-12}$ \\ 
 2 &   1.2830$\times$10$^{4}$ &   1.2830$\times$10$^{4}$ &  1.0000$\times$10$^{-16}$ \\ 
 3 &   2.5152$\times$10$^{4}$ &   2.5152$\times$10$^{4}$ &  7.2760$\times$10$^{-12}$ \\ 
 4 &   4.1578$\times$10$^{4}$ &   4.1578$\times$10$^{4}$ &  1.0000$\times$10$^{-16}$ \\ 
 5 &   6.2110$\times$10$^{4}$ &   6.2110$\times$10$^{4}$ &  1.8190$\times$10$^{-10}$ \\
\bottomrule
\end{tabular}
\end{threeparttable}
\end{table}\\
\newline
As seen in Table~\ref{tab:1d_transverse_fixed}, the eigenfrequencies computed using Eig-PIELM closely match the analytical values for all five modes considered. The absolute differences are extremely small, ranging from $1.0\times 10^{-16}$ to $1.81\times 10^{-10},\mathrm{rad/s}$, which correspond to relative errors in the range of $10^{-16}$ to $10^{-13}$. Such small discrepancies show that the method reproduces both the fundamental and higher transverse modes with very high accuracy in this setting. Using $N_\phi = 24$ neurons and $N_x = 50{,}000$ interior points, the full set of modes was obtained in $0.16,\mathrm{s}$.
The consistently low errors across modes suggest that Eig-PIELM maintains its reliability even for higher-order, more oscillatory mode shapes in fixed–fixed beams.

\paragraph{Fixed-free} As a last case, consider a homogeneous Dirichlet and zero-slope condition are imposed at the fixed end $x=0$, while the free end $x=L$ satisfies zero bending moment and zero shear force, i.e.,
\begin{equation}
\left\{ 
\begin{aligned}
    &w(0) = 0, ~~~~~~\qquad \dfrac{dw}{dx}\Big|_{x=0} = 0, \\
    &\dfrac{d^2w}{dx^2}\Big|_{x=L} = 0, \quad \dfrac{d^3w}{dx^3}\Big|_{x=L} = 0.
\end{aligned}
\right.
\end{equation}
The exact natural frequencies \cite{rao2007vibration} are
\begin{equation}
    \omega_n = \left( \frac{\beta_n}{L} \right)^2
    \sqrt{\frac{E I}{\rho A}}, 
    \quad n = 1, 2, \dots,
\end{equation}
where $\beta_n = \{1.8751,\, 4.6941,\, 7.8548,\, 10.9955,\, 14.1372\}$ are the dimensionless wavenumbers for the first five cantilever modes, determined from the characteristic equation. The corresponding normalized mode shapes are given \cite{rao2007vibration} by
\begin{equation}
    \phi_n(x) = \left[ \cos\!\left( \frac{\beta_n x}{L} \right) - \cosh\!\left( \frac{\beta_n x}{L} \right)\right]
    - \left( \dfrac{\cos\beta_n + \cosh\beta_n}{\sin\beta_n + \sinh\beta_n} \right) \left[ \sin\!\left( \frac{\beta_n x}{L} \right) -\sinh\!\left( \frac{\beta_n x}{L} \right) \right],
\end{equation}
\Cref{fig:1d_transverse_fixed_free} shows the first five mode shapes, comparing Eig-PIELM predictions with the exact analytical solutions. Table~\ref{tab:1d_transverse_cantilever} presents the corresponding eigenfrequencies.
\begin{table}[!ht]
\centering
\normalsize
\begin{threeparttable}
\captionsetup{labelformat=empty, font=small}
\caption{\textbf{Table \thetable} \\ Comparison of predicted and analytical eigenfrequencies for the first five transverse vibration modes of a one-dimensional Euler–Bernoulli cantilever beam (fixed–free)}
\label{tab:1d_transverse_cantilever}
\begin{tabular}{lc@{\hspace{2em}}c@{\hspace{2em}}c@{\hspace{2em}}c@{\hspace{2em}}c}
\toprule
Mode & $\omega_{\scriptscriptstyle\text{NP-PIELM}}$ & $\omega_{\scriptscriptstyle\text{Exact}}$ & $|\omega_{\scriptscriptstyle\text{NP-PIELM}}-\omega_{\scriptscriptstyle\text{Exact}}|$ \\
Number & $\scriptstyle(\text{rad/s})$ & $\scriptstyle(\text{rad/s})$ & $\scriptstyle(\text{rad/s})$ \\
\midrule
 1 &   7.3146$\times$10$^{2}$ &   7.3146$\times$10$^{2}$ &  2.1687$\times$10$^{-09}$ \\ 
 2 &   4.5840$\times$10$^{3}$ &   4.5840$\times$10$^{3}$ &  1.5280$\times$10$^{-10}$ \\ 
 3 &   1.2835$\times$10$^{4}$ &   1.2835$\times$10$^{4}$ &  3.0195$\times$10$^{-10}$ \\ 
 4 &   2.5152$\times$10$^{4}$ &   2.5152$\times$10$^{4}$ &  1.6335$\times$10$^{-09}$ \\ 
 5 &   4.1578$\times$10$^{4s}$ &   4.1578$\times$10$^{4}$ &  8.8374$\times$10$^{-08}$ \\
\bottomrule
\end{tabular}
\end{threeparttable}
\end{table}
\noindent
As shown in Table~\ref{tab:1d_transverse_cantilever}, the predicted eigenfrequencies are in extremely close agreement with the analytical values. Absolute differences range from $1.53\times 10^{-10}$ to $8.84\times 10^{-8}\,\mathrm{rad/s}$, corresponding to relative errors of approximately $10^{-14}$ to $10^{-12}$. Both the fundamental and higher-order modes are reproduced with very high accuracy. Using $N_{\phi} = 24$ neurons and $N_x = 50{,}000$ interior points, the computation completed in $0.177\,\mathrm{s}$. The largest deviation occurs for the fifth mode, which is expected due to its higher oscillatory content, but the accuracy remains excellent across all modes.


\subsection{Two dimensional rigid-walled rectangular acoustic cavity}
\begin{figure}
    \centering
    \includegraphics[width=0.9\linewidth,trim=150 20 60 20,clip]{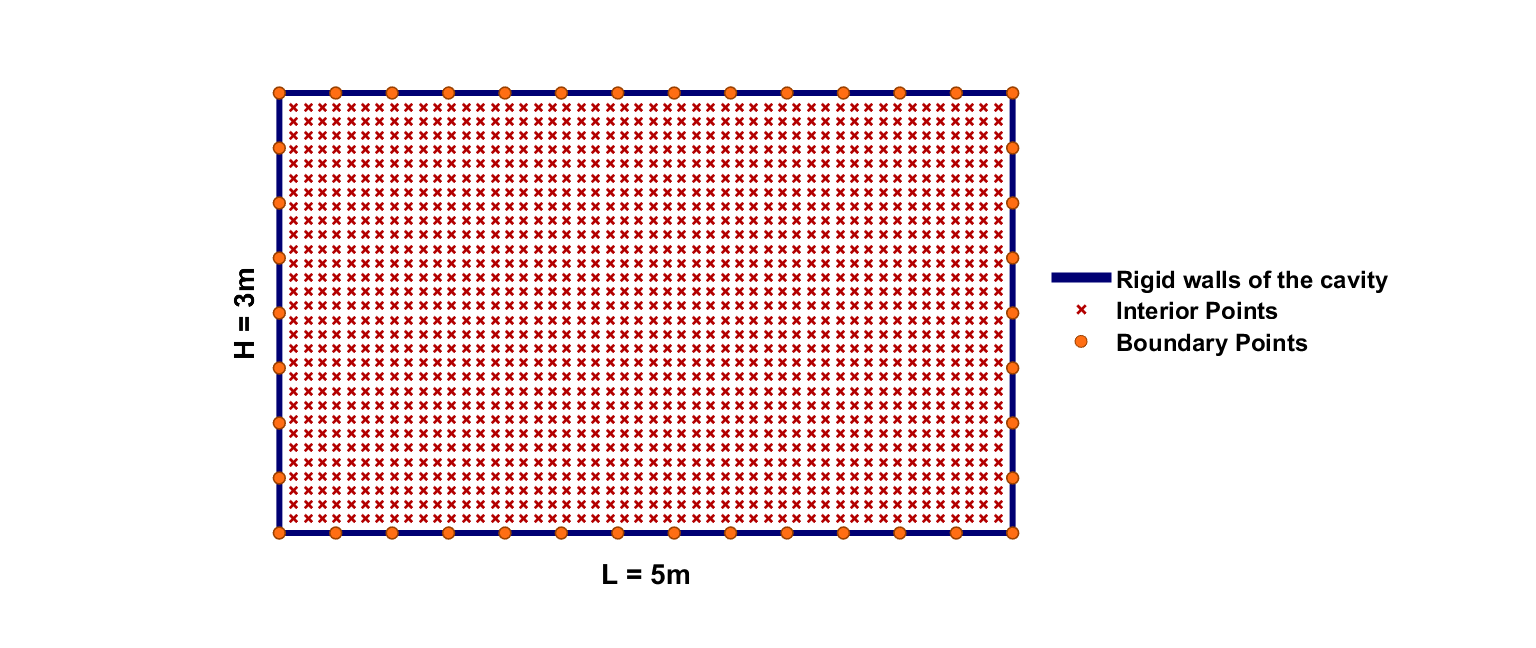}
    \caption{Distribution of interior and boundary sampling points for the 2D rigid-walled rectangular cavity.}
    \label{fig:2D_cavity_domain_points}
    
\end{figure}

\begin{figure}
    \centering
    \includegraphics[width=1.1\linewidth]{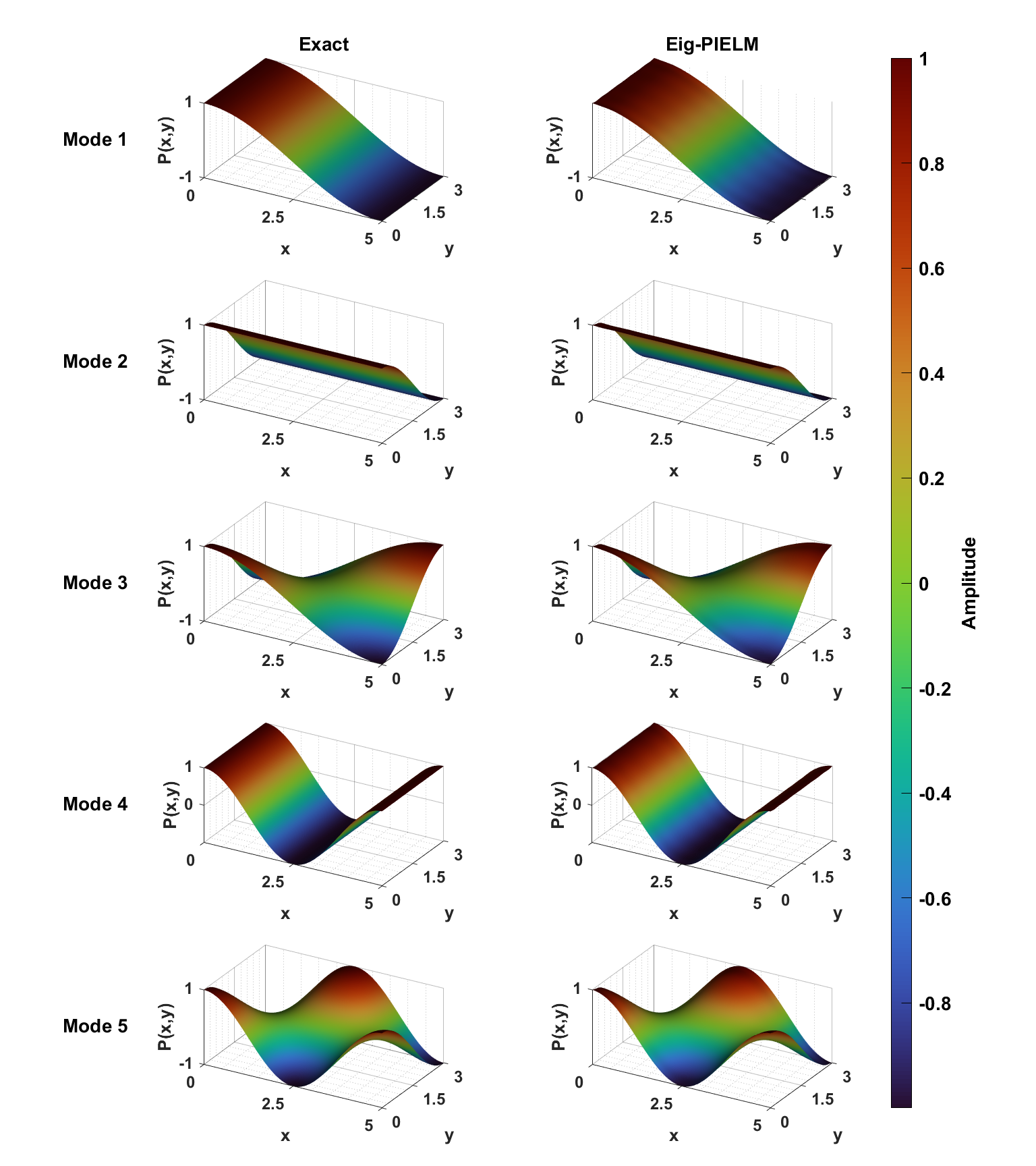}
    \caption{First five mode shapes of a 2D rigid-walled rectangular cavity from Eig-PIELM and exact analytical solutions, shown as 3D surface plots of the acoustic pressure field $P(x,y)$.}

    \label{fig:2D_Helmhotz_mode_shapes}
\end{figure}

\noindent In this example, we solve the homogeneous Helmholtz equation, see \Cref{eqn:helm_freq} on a rectangular domain, $\Omega = [0,L] \times [0,H]$ with Neumann conditions enforced on the boundary. The analytical eigenpairs for a rectangular domain are doubly indexed:
\begin{equation}
\begin{aligned}
    P_{m,n}(x,y) &= \cos\!\left( \frac{m\pi x}{L} \right)\cos\!\left( \frac{n\pi y}{H} \right),\\
    k_{m,n}^2 &= \left( \frac{m\pi}{L} \right)^2 + \left( \frac{n\pi}{H} \right)^2, 
    \quad m,n = 0,1,2,\dots
\end{aligned}
\end{equation}
with $(m,n) \neq (0,0)$. The corresponding natural frequencies are $\omega_{m,n} = c\,k_{m,n}$. In the Eig-PIELM implementation, $P(x,y)$ is approximated using a polynomial basis of degree $n_x = 13$ in $x$ and $n_y= 13$ in $y$ directions, resulting in a total of $(n_x+1)(n_y+1) = 196$ basis functions. The governing PDE (c.f. \Cref{eqn:helm_freq} is enforced at $N_x = 500$ sample points in the $x$-direction and $N_y = 3000$ points in the $y$-direction. The complete eigenspectrum is computed in $0.244\,\mathrm{s}$. \Cref{fig:2D_cavity_domain_points} shows the distribution of interior and boundary sample points. \Cref{fig:2D_Helmhotz_mode_shapes} presents the first five mode shapes as side-by-side 3D surface plots of $P(x,y)$, displaying the Eig-PIELM predictions alongside the exact solutions.
\begin{table}[!ht]
\centering
\normalsize
\begin{threeparttable}
\captionsetup{labelformat=empty, font=small}
\caption{\textbf{Table \thetable} \\Comparison of predicted and exact eigenfrequencies (rad/s) for the first five modes of a 2D rigid-walled rectangular cavity.}
\label{tab:2d_helmhotz}
\begin{tabular}{lc@{\hspace{2em}}c@{\hspace{2em}}c@{\hspace{2em}}c@{\hspace{2em}}c}
\toprule
Mode & $\omega_{\scriptscriptstyle\text{NP-PIELM}}$ & $\omega_{\scriptscriptstyle\text{Exact}}$ & $|\omega_{\scriptscriptstyle\text{NP-PIELM}}-\omega_{\scriptscriptstyle\text{Exact}}|$ \\
Number & $\scriptstyle(\text{rad/s})$ & $\scriptstyle(\text{rad/s})$ & $\scriptstyle(\text{rad/s})$ \\
\midrule
 1 &  2.1363$\times$10$^{2}$ &   2.1363$\times$10$^{2}$ &  6.1890$\times$10$^{-05}$ \\ 
 2 &  3.5605$\times$10$^{2}$ &   3.5605$\times$10$^{2}$ &  1.2857$\times$10$^{-05}$ \\ 
 3 &  4.1522$\times$10$^{2}$ &   4.1522$\times$10$^{2}$ &  9.0340$\times$10$^{-05}$ \\ 
 4 &  4.2726$\times$10$^{2}$ &   4.2726$\times$10$^{2}$ &  6.625e$\times$10$^{-06}$ \\ 
 5 &  5.5616$\times$10$^{2}$ &   5.5616$\times$10$^{2}$ &  9.3155$\times$10$^{-06}$ \\  
\bottomrule
\end{tabular}
\end{threeparttable}
\end{table}
\newline\\
Table~\ref{tab:2d_helmhotz} indicates close agreement between Eig-PIELM and the exact eigenfrequencies. Absolute discrepancies for the first five modes range from $6.63\times 10^{-6}$ to $9.03\times 10^{-5}\,\text{rad/s}$. In relative terms, the errors are in the range of \(\mathcal{O}(10^{-8})\) \text{ to } \(\mathcal{O}(10^{-7})\). Overall, the results show that Eig-PIELM reproduces both lower- and higher-index modes with very high accuracy at modest computational cost.

\section{Conclusions}\label{sec:conc}
\noindent
This work successfully extends the Physics-Informed Extreme Learning Machine (PIELM) framework to solve linear eigenvalue problems through the novel Eig-PIELM approach. The fundamental challenge addressed was that standard PIELM struggles with eigenvalue problems because the eigenvalue and the eigenfunctions are unknown simultaneously, which renders the stationarity conditions nonlinear. This nonlinearity disrupts the closed-form solution paradigm that defines PIELM's computational advantage.\\
\newline
The proposed Eig-PIELM framework resolves this  issue by enforcing boundary conditions exactly in coefficient space and collapsing the two stationarity equations to a single symmetric generalized eigen value problem. This strategic formulation preserves the usual strengths of PIELM including mesh-free sampling, elimination of backpropagation, and direct linear algebra while also enabling entire spectrum of modes to be computed in one solve, rather than requiring separate optimization runs for each eigenmode.\\
\newline
Comprehensive numerical validation across diverse benchmark problems, including transverse vibrations of one-dimensional beams and acoustic cavity eigenanalysis, demonstrates exceptional accuracy. For one-dimensional free vibration of beams, the relative errors  are at machine precision $\sim10^{-16}$ for fundamental modes to better than $10^{-9}$ for higher-order modes. The two-dimensional Helmholtz equation on a rectangular acoustic cavity produces relative errors of $10^{-8}$ to $10^{-7}$. The computational efficiency is remarkable, with complete eigenspectra obtained in 0.16 to 0.244 seconds, representing potential
speedups over iterative methods.\\
\newline 
The exact boundary condition enforcement, combined with the emergent symmetrization property in the boundary-admissible coordinate space, establishes Eig-PIELM as a robust and theoretically sound framework.\\
\newline
In our opinion, Eig-PIELM fills a critical gap in mesh-free eigenanalysis, particularly for parametric studies and design optimization requiring rapid, accurate eigenvalue predictions. Future extensions to nonlinear eigenvalue problems and large-scale implementations will broaden its applicability to industrial applications.

\bibliographystyle{elsarticle-num-names} 
\bibliography{myRef}

\begin{thebibliography}{55}
\expandafter\ifx\csname natexlab\endcsname\relax\def\natexlab#1{#1}\fi
\providecommand{\url}[1]{\texttt{#1}}
\providecommand{\href}[2]{#2}
\providecommand{\path}[1]{#1}
\providecommand{\DOIprefix}{doi:}
\providecommand{\ArXivprefix}{arXiv:}
\providecommand{\URLprefix}{URL: }
\providecommand{\Pubmedprefix}{pmid:}
\providecommand{\doi}[1]{\href{http://dx.doi.org/#1}{\path{#1}}}
\providecommand{\Pubmed}[1]{\href{pmid:#1}{\path{#1}}}
\providecommand{\bibinfo}[2]{#2}
\ifx\xfnm\relax \def\xfnm[#1]{\unskip,\space#1}\fi
\bibitem[{Schr{\"o}dinger(1940)}]{schrodinger1940method}
\bibinfo{author}{E.~Schr{\"o}dinger},
\newblock \bibinfo{title}{A method of determining quantum-mechanical eigenvalues and eigenfunctions},
\newblock in: \bibinfo{booktitle}{Proceedings of the Royal Irish Academy. Section A: Mathematical and Physical Sciences}, volume~\bibinfo{volume}{46}, \bibinfo{organization}{JSTOR}, \bibinfo{year}{1940}, pp. \bibinfo{pages}{9--16}.
\bibitem[{Turbiner(1984)}]{turbiner1984eigenvalue}
\bibinfo{author}{A.~Turbiner},
\newblock \bibinfo{title}{The eigenvalue spectrum in quantum mechanics and the nonlinearizationprocedure},
\newblock \bibinfo{journal}{Soviet Physics Uspekhi} \bibinfo{volume}{27} (\bibinfo{year}{1984}) \bibinfo{pages}{668}.
\bibitem[{Bathe and Wilson(1973)}]{bathe1973solution}
\bibinfo{author}{K.-J. Bathe}, \bibinfo{author}{E.~L. Wilson},
\newblock \bibinfo{title}{Solution methods for eigenvalue problems in structural mechanics},
\newblock \bibinfo{journal}{International Journal for Numerical Methods in Engineering} \bibinfo{volume}{6} (\bibinfo{year}{1973}) \bibinfo{pages}{213--226}.
\bibitem[{Shinozuka and Astill(1972)}]{shinozuka1972random}
\bibinfo{author}{M.~Shinozuka}, \bibinfo{author}{C.~J. Astill},
\newblock \bibinfo{title}{Random eigenvalue problems in structural analysis},
\newblock \bibinfo{journal}{AIAA journal} \bibinfo{volume}{10} (\bibinfo{year}{1972}) \bibinfo{pages}{456--462}.
\bibitem[{Cliffe et~al.(1994)Cliffe, Garratt, and Spence}]{cliffe1994eigenvalues}
\bibinfo{author}{K.~A. Cliffe}, \bibinfo{author}{T.~J. Garratt}, \bibinfo{author}{A.~Spence},
\newblock \bibinfo{title}{Eigenvalues of block matrices arising from problems in fluid mechanics},
\newblock \bibinfo{journal}{SIAM Journal on Matrix Analysis and Applications} \bibinfo{volume}{15} (\bibinfo{year}{1994}) \bibinfo{pages}{1310--1318}.
\bibitem[{Mahajan et~al.(1991)Mahajan, Dowell, and Bliss}]{mahajan1991eigenvalue}
\bibinfo{author}{A.~J. Mahajan}, \bibinfo{author}{E.~H. Dowell}, \bibinfo{author}{D.~B. Bliss},
\newblock \bibinfo{title}{Eigenvalue calculation procedure for an euler/navier-stokes solver with application to flows over airfoils},
\newblock \bibinfo{journal}{Journal of Computational Physics} \bibinfo{volume}{97} (\bibinfo{year}{1991}) \bibinfo{pages}{398--413}.
\bibitem[{Valovik(2018)}]{valovik2018nonlinear}
\bibinfo{author}{D.~Valovik},
\newblock \bibinfo{title}{On a nonlinear eigenvalue problem related to the theory of propagation of electromagnetic waves},
\newblock \bibinfo{journal}{Differential Equations} \bibinfo{volume}{54} (\bibinfo{year}{2018}) \bibinfo{pages}{165--177}.
\bibitem[{Guo(2009)}]{guo2009eigen}
\bibinfo{author}{S.~Guo},
\newblock \bibinfo{title}{An eigen theory of electromagnetic waves based on the standard spaces},
\newblock \bibinfo{journal}{International Journal of Engineering Science} \bibinfo{volume}{47} (\bibinfo{year}{2009}) \bibinfo{pages}{405--412}.
\bibitem[{Strang(1993)}]{strang1993linear}
\bibinfo{author}{G.~Strang}, \bibinfo{title}{Introduction to Linear Algebra}, \bibinfo{edition}{1} ed., \bibinfo{publisher}{Wellesley--Cambridge Press}, \bibinfo{year}{1993}.
\bibitem[{Horn and Johnson(2013)}]{horn2013matrix}
\bibinfo{author}{R.~A. Horn}, \bibinfo{author}{C.~R. Johnson}, \bibinfo{title}{Matrix Analysis}, \bibinfo{edition}{2} ed., \bibinfo{publisher}{Cambridge University Press}, \bibinfo{year}{2013}.
\bibitem[{Choi and Hoefer(1986)}]{choi1986finite}
\bibinfo{author}{D.~H. Choi}, \bibinfo{author}{W.~J. Hoefer},
\newblock \bibinfo{title}{The finite-difference-time-domain method and its application to eigenvalue problems},
\newblock \bibinfo{journal}{IEEE Transactions on Microwave Theory and Techniques} \bibinfo{volume}{34} (\bibinfo{year}{1986}) \bibinfo{pages}{1464--1470}.
\bibitem[{Baxley(1972)}]{baxley1972eigenvalues}
\bibinfo{author}{J.~V. Baxley},
\newblock \bibinfo{title}{Eigenvalues of singular differential operators by finite difference methods, i},
\newblock \bibinfo{journal}{Journal of Mathematical Analysis and Applications} \bibinfo{volume}{38} (\bibinfo{year}{1972}) \bibinfo{pages}{244--254}.
\bibitem[{Bertrand et~al.(2023)Bertrand, Boffi, and Halim}]{bertrand2023reduced}
\bibinfo{author}{F.~Bertrand}, \bibinfo{author}{D.~Boffi}, \bibinfo{author}{A.~Halim},
\newblock \bibinfo{title}{A reduced order model for the finite element approximation of eigenvalue problems},
\newblock \bibinfo{journal}{Computer Methods in Applied Mechanics and Engineering} \bibinfo{volume}{404} (\bibinfo{year}{2023}) \bibinfo{pages}{115696}.
\bibitem[{Hughes et~al.(2014)Hughes, Evans, and Reali}]{hughes2014finite}
\bibinfo{author}{T.~J. Hughes}, \bibinfo{author}{J.~A. Evans}, \bibinfo{author}{A.~Reali},
\newblock \bibinfo{title}{Finite element and nurbs approximations of eigenvalue, boundary-value, and initial-value problems},
\newblock \bibinfo{journal}{Computer Methods in Applied Mechanics and Engineering} \bibinfo{volume}{272} (\bibinfo{year}{2014}) \bibinfo{pages}{290--320}.
\bibitem[{Chen et~al.(2011)Chen, He, and Zhou}]{chen2011finite}
\bibinfo{author}{H.~Chen}, \bibinfo{author}{L.~He}, \bibinfo{author}{A.~Zhou},
\newblock \bibinfo{title}{Finite element approximations of nonlinear eigenvalue problems in quantum physics},
\newblock \bibinfo{journal}{Computer methods in applied mechanics and engineering} \bibinfo{volume}{200} (\bibinfo{year}{2011}) \bibinfo{pages}{1846--1865}.
\bibitem[{Feit et~al.(1982)Feit, Fleck~Jr, and Steiger}]{feit1982solution}
\bibinfo{author}{M.~Feit}, \bibinfo{author}{J.~Fleck~Jr}, \bibinfo{author}{A.~Steiger},
\newblock \bibinfo{title}{Solution of the schr{\"o}dinger equation by a spectral method},
\newblock \bibinfo{journal}{Journal of Computational Physics} \bibinfo{volume}{47} (\bibinfo{year}{1982}) \bibinfo{pages}{412--433}.
\bibitem[{Gheorghiu(2014)}]{gheorghiu2014spectral}
\bibinfo{author}{C.-I. Gheorghiu}, \bibinfo{title}{Spectral methods for non-standard eigenvalue problems: fluid and structural mechanics and beyond}, \bibinfo{publisher}{Springer Science \& Business}, \bibinfo{year}{2014}.
\bibitem[{Canuto et~al.(2007)Canuto, Hussaini, Quarteroni, and Zang}]{canuto2007spectral}
\bibinfo{author}{C.~Canuto}, \bibinfo{author}{M.~Y. Hussaini}, \bibinfo{author}{A.~Quarteroni}, \bibinfo{author}{T.~A. Zang}, \bibinfo{title}{{Spectral Methods: Fundamentals in Single Domains}}, \bibinfo{publisher}{Springer}, \bibinfo{year}{2007}.
\bibitem[{Lanczos(1950)}]{Lanczos1950}
\bibinfo{author}{C.~Lanczos},
\newblock \bibinfo{title}{An iteration method for the solution of the eigenvalue problem of linear differential and integral operators},
\newblock \bibinfo{journal}{Journal of Research of the National Bureau of Standards} \bibinfo{volume}{45} (\bibinfo{year}{1950}) \bibinfo{pages}{255--282}.
\bibitem[{Francis(1961)}]{Francis1961}
\bibinfo{author}{J.~G.~F. Francis},
\newblock \bibinfo{title}{The qr transformation, part ii},
\newblock \bibinfo{journal}{Computer Journal} \bibinfo{volume}{4} (\bibinfo{year}{1961}) \bibinfo{pages}{332--345}.
\bibitem[{Parlett(1998)}]{parlett1998symmetric}
\bibinfo{author}{B.~Parlett}, \bibinfo{title}{The Symmetric Eigenvalue Problem}, \bibinfo{publisher}{SIAM}, \bibinfo{year}{1998}.
\bibitem[{Trefethen and Bau(1997)}]{trefethen1997numerical}
\bibinfo{author}{L.~Trefethen}, \bibinfo{author}{D.~Bau}, \bibinfo{title}{Numerical Linear Algebra}, \bibinfo{publisher}{SIAM}, \bibinfo{year}{1997}.
\bibitem[{Saad(1980)}]{Saad1980}
\bibinfo{author}{Y.~Saad},
\newblock \bibinfo{title}{Variations on arnoldi’s method for computing eigenvalues of large unsymmetric matrices},
\newblock \bibinfo{journal}{Linear Algebra and Its Applications} \bibinfo{volume}{34} (\bibinfo{year}{1980}) \bibinfo{pages}{269--295}.
\bibitem[{Lehoucq et~al.(1998)Lehoucq, Sorensen, and Yang}]{Lehoucq1998}
\bibinfo{author}{R.~B. Lehoucq}, \bibinfo{author}{D.~C. Sorensen}, \bibinfo{author}{C.~Yang}, \bibinfo{title}{ARPACK Users’ Guide: Solution of Large‐Scale Eigenvalue Problems with Implicitly Restarted Arnoldi Methods}, \bibinfo{publisher}{SIAM}, \bibinfo{year}{1998}.
\bibitem[{Saad(2011)}]{saad2011numerical}
\bibinfo{author}{Y.~Saad}, \bibinfo{title}{Numerical Methods for Large Eigenvalue Problems}, \bibinfo{edition}{2} ed., \bibinfo{publisher}{SIAM}, \bibinfo{year}{2011}.
\bibitem[{Lagaris et~al.(1997)Lagaris, Likas, and Fotiadis}]{lagaris1997artificial}
\bibinfo{author}{I.~E. Lagaris}, \bibinfo{author}{A.~Likas}, \bibinfo{author}{D.~I. Fotiadis},
\newblock \bibinfo{title}{Artificial neural network methods in quantum mechanics},
\newblock \bibinfo{journal}{Computer Physics Communications} \bibinfo{volume}{104} (\bibinfo{year}{1997}) \bibinfo{pages}{1--14}.
\bibitem[{Yu et~al.(2018)}]{yu2018deep}
\bibinfo{author}{B.~Yu}, et~al.,
\newblock \bibinfo{title}{The deep ritz method: a deep learning-based numerical algorithm for solving variational problems},
\newblock \bibinfo{journal}{Communications in Mathematics and Statistics} \bibinfo{volume}{6} (\bibinfo{year}{2018}) \bibinfo{pages}{1--12}.
\bibitem[{Raissi et~al.(2019)Raissi, Perdikaris, and Karniadakis}]{raissi2019physics}
\bibinfo{author}{M.~Raissi}, \bibinfo{author}{P.~Perdikaris}, \bibinfo{author}{G.~E. Karniadakis},
\newblock \bibinfo{title}{Physics-informed neural networks: A deep learning framework for solving forward and inverse problems involving nonlinear partial differential equations},
\newblock \bibinfo{journal}{Journal of Computational physics} \bibinfo{volume}{378} (\bibinfo{year}{2019}) \bibinfo{pages}{686--707}.
\bibitem[{Jin et~al.(2020)Jin, Mattheakis, and Protopapas}]{jin2020unsupervised}
\bibinfo{author}{H.~Jin}, \bibinfo{author}{M.~Mattheakis}, \bibinfo{author}{P.~Protopapas},
\newblock \bibinfo{title}{Unsupervised neural networks for quantum eigenvalue problems},
\newblock \bibinfo{journal}{arXiv preprint arXiv:2010.05075}  (\bibinfo{year}{2020}).
\bibitem[{Jin et~al.(2022)Jin, Mattheakis, and Protopapas}]{jin2022physics}
\bibinfo{author}{H.~Jin}, \bibinfo{author}{M.~Mattheakis}, \bibinfo{author}{P.~Protopapas},
\newblock \bibinfo{title}{Physics-informed neural networks for quantum eigenvalue problems},
\newblock in: \bibinfo{booktitle}{2022 International Joint Conference on Neural Networks (IJCNN)}, \bibinfo{organization}{IEEE}, \bibinfo{year}{2022}, pp. \bibinfo{pages}{1--8}.
\bibitem[{Ben-Shaul et~al.(2023)Ben-Shaul, Bar, Fishelov, and Sochen}]{ben2023deep}
\bibinfo{author}{I.~Ben-Shaul}, \bibinfo{author}{L.~Bar}, \bibinfo{author}{D.~Fishelov}, \bibinfo{author}{N.~Sochen},
\newblock \bibinfo{title}{Deep learning solution of the eigenvalue problem for differential operators},
\newblock \bibinfo{journal}{Neural Computation} \bibinfo{volume}{35} (\bibinfo{year}{2023}) \bibinfo{pages}{1100--1134}.
\bibitem[{Yang et~al.(2023)Yang, Gong, Zhang, Yang, Chen, He, and Li}]{yang2023data}
\bibinfo{author}{Y.~Yang}, \bibinfo{author}{H.~Gong}, \bibinfo{author}{S.~Zhang}, \bibinfo{author}{Q.~Yang}, \bibinfo{author}{Z.~Chen}, \bibinfo{author}{Q.~He}, \bibinfo{author}{Q.~Li},
\newblock \bibinfo{title}{A data-enabled physics-informed neural network with comprehensive numerical study on solving neutron diffusion eigenvalue problems},
\newblock \bibinfo{journal}{Annals of Nuclear Energy} \bibinfo{volume}{183} (\bibinfo{year}{2023}) \bibinfo{pages}{109656}.
\bibitem[{Yu et~al.(2024)Yu, He, Zhang, Yang, Yang, and Gong}]{yu2024solving}
\bibinfo{author}{S.~Yu}, \bibinfo{author}{Q.~He}, \bibinfo{author}{S.~Zhang}, \bibinfo{author}{Q.~Yang}, \bibinfo{author}{Y.~Yang}, \bibinfo{author}{H.~Gong},
\newblock \bibinfo{title}{Solving multi-group neutron diffusion eigenvalue problem with decoupling residual loss function},
\newblock \bibinfo{journal}{arXiv preprint arXiv:2411.15693}  (\bibinfo{year}{2024}).
\bibitem[{Zhang et~al.(2024)Zhang, Jiang, Qian, and Xu}]{zhang2024orthogonal}
\bibinfo{author}{S.~Zhang}, \bibinfo{author}{X.~Jiang}, \bibinfo{author}{H.~Qian}, \bibinfo{author}{Y.~Xu},
\newblock \bibinfo{title}{Orthogonal constrained neural networks for solving structured inverse eigenvalue problems},
\newblock \bibinfo{journal}{arXiv preprint arXiv:2406.19981}  (\bibinfo{year}{2024}).
\bibitem[{Pallikarakis and Ntargaras(2024)}]{pallikarakis2024application}
\bibinfo{author}{N.~Pallikarakis}, \bibinfo{author}{A.~Ntargaras},
\newblock \bibinfo{title}{Application of machine learning regression models to inverse eigenvalue problems},
\newblock \bibinfo{journal}{Computers \& Mathematics with Applications} \bibinfo{volume}{154} (\bibinfo{year}{2024}) \bibinfo{pages}{162--174}.
\bibitem[{Singhal and Agarwal(2024)}]{singhal2024physics}
\bibinfo{author}{A.~Singhal}, \bibinfo{author}{H.~Agarwal},
\newblock \bibinfo{title}{Physics informed neural network based time-independent schr{\"o}dinger equation solver},
\newblock in: \bibinfo{booktitle}{2024 8th IEEE Electron Devices Technology \& Manufacturing Conference (EDTM)}, \bibinfo{organization}{IEEE}, \bibinfo{year}{2024}, pp. \bibinfo{pages}{1--3}.
\bibitem[{Mattheakis et~al.(2022)Mattheakis, Schleder, Larson, and Kaxiras}]{mattheakis2022first}
\bibinfo{author}{M.~Mattheakis}, \bibinfo{author}{G.~R. Schleder}, \bibinfo{author}{D.~T. Larson}, \bibinfo{author}{E.~Kaxiras},
\newblock \bibinfo{title}{First principles physics-informed neural network for quantum wavefunctions and eigenvalue surfaces},
\newblock \bibinfo{journal}{arXiv preprint arXiv:2211.04607}  (\bibinfo{year}{2022}). \href{http://arxiv.org/abs/2211.04607}{{\tt arXiv:2211.04607}}.
\bibitem[{Chen et~al.(2024)Chen, Lai, and Yang}]{chen2024pinn}
\bibinfo{author}{Z.~Chen}, \bibinfo{author}{S.-K. Lai}, \bibinfo{author}{Z.~Yang},
\newblock \bibinfo{title}{At-pinn: Advanced time-marching physics-informed neural network for structural vibration analysis},
\newblock \bibinfo{journal}{Thin-Walled Structures} \bibinfo{volume}{196} (\bibinfo{year}{2024}) \bibinfo{pages}{111423}.
\bibitem[{Yoo et~al.(2025)Yoo, Kang, Yoon, and Kim}]{yoo2025physics}
\bibinfo{author}{S.~Yoo}, \bibinfo{author}{M.~Kang}, \bibinfo{author}{H.~Yoon}, \bibinfo{author}{T.~Kim},
\newblock \bibinfo{title}{A physics-informed neural network approach for solving structural eigenvalue problem},
\newblock \bibinfo{journal}{International Journal of Precision Engineering and Manufacturing}  (\bibinfo{year}{2025}) \bibinfo{pages}{1--14}.
\bibitem[{Krishnapriyan et~al.(2021)Krishnapriyan, Gholami, Zhe, Kirby, and Mahoney}]{krishnapriyan2021characterizing}
\bibinfo{author}{A.~Krishnapriyan}, \bibinfo{author}{A.~Gholami}, \bibinfo{author}{S.~Zhe}, \bibinfo{author}{R.~Kirby}, \bibinfo{author}{M.~W. Mahoney},
\newblock \bibinfo{title}{Characterizing possible failure modes in physics-informed neural networks},
\newblock \bibinfo{journal}{Advances in Neural Information Processing Systems} \bibinfo{volume}{34} (\bibinfo{year}{2021}) \bibinfo{pages}{26548--26560}.
\bibitem[{Wang et~al.(2021)Wang, Teng, and Perdikaris}]{wang2022understanding}
\bibinfo{author}{S.~Wang}, \bibinfo{author}{Y.~Teng}, \bibinfo{author}{P.~Perdikaris},
\newblock \bibinfo{title}{Understanding and mitigating gradient flow pathologies in physics-informed neural networks},
\newblock \bibinfo{journal}{SIAM Journal on Scientific Computing} \bibinfo{volume}{43} (\bibinfo{year}{2021}) \bibinfo{pages}{A3055--A3081}.
\bibitem[{Grossmann et~al.(2024)Grossmann, Komorowska, Latz, and Sch{\"o}nlieb}]{grossmann2023can}
\bibinfo{author}{T.~G. Grossmann}, \bibinfo{author}{U.~J. Komorowska}, \bibinfo{author}{J.~Latz}, \bibinfo{author}{C.-B. Sch{\"o}nlieb},
\newblock \bibinfo{title}{Can physics-informed neural networks beat the finite element method?},
\newblock \bibinfo{journal}{IMA Journal of Applied Mathematics} \bibinfo{volume}{89} (\bibinfo{year}{2024}) \bibinfo{pages}{143--187}.
\bibitem[{Rahaman et~al.(2019)Rahaman, Baratin, Arpit, Draxler, Lin, Hamprecht, Bengio, and Courville}]{rahaman2019spectral}
\bibinfo{author}{N.~Rahaman}, \bibinfo{author}{A.~Baratin}, \bibinfo{author}{D.~Arpit}, \bibinfo{author}{F.~Draxler}, \bibinfo{author}{M.~Lin}, \bibinfo{author}{F.~A. Hamprecht}, \bibinfo{author}{Y.~Bengio}, \bibinfo{author}{A.~Courville},
\newblock \bibinfo{title}{On the spectral bias of neural networks},
\newblock \bibinfo{journal}{International conference on machine learning}  (\bibinfo{year}{2019}) \bibinfo{pages}{5301--5310}.
\bibitem[{Tancik et~al.(2020)Tancik, Srinivasan, Mildenhall, Fridovich-Keil, Raghavan, Singhal, Ramamoorthi, Barron, and Ng}]{tancik2020fourier}
\bibinfo{author}{M.~Tancik}, \bibinfo{author}{P.~P. Srinivasan}, \bibinfo{author}{B.~Mildenhall}, \bibinfo{author}{S.~Fridovich-Keil}, \bibinfo{author}{N.~Raghavan}, \bibinfo{author}{U.~Singhal}, \bibinfo{author}{R.~Ramamoorthi}, \bibinfo{author}{J.~T. Barron}, \bibinfo{author}{R.~Ng},
\newblock \bibinfo{title}{Fourier features let networks learn high frequency functions in low dimensional domains},
\newblock \bibinfo{journal}{Advances in Neural Information Processing Systems} \bibinfo{volume}{33} (\bibinfo{year}{2020}) \bibinfo{pages}{7537--7547}.
\bibitem[{Wang et~al.(2021)Wang, Hanwen, and Perdikaris}]{wang2021eigenvector}
\bibinfo{author}{S.~Wang}, \bibinfo{author}{W.~Hanwen}, \bibinfo{author}{P.~Perdikaris},
\newblock \bibinfo{title}{On the eigenvector bias of fourier feature networks: From regression to solving multi-scale pdes with physics-informed neural networks},
\newblock \bibinfo{journal}{Computer Methods in Applied Mechanics and Engineering} \bibinfo{volume}{384} (\bibinfo{year}{2021}) \bibinfo{pages}{113938}.
\bibitem[{Dwivedi and Srinivasan(2020)}]{dwivedi2020physics}
\bibinfo{author}{V.~Dwivedi}, \bibinfo{author}{B.~Srinivasan},
\newblock \bibinfo{title}{Physics informed extreme learning machine (pielm)--a rapid method for the numerical solution of partial differential equations},
\newblock \bibinfo{journal}{Neurocomputing} \bibinfo{volume}{391} (\bibinfo{year}{2020}) \bibinfo{pages}{96--118}.
\bibitem[{Pan and Li(2024)}]{pan2024applying}
\bibinfo{author}{S.~Pan}, \bibinfo{author}{H.~Li},
\newblock \bibinfo{title}{Applying physics-informed extreme learning machines to solve the euler-bernoulli beam problem},
\newblock \bibinfo{journal}{Theoretical and Natural Science} \bibinfo{volume}{53} (\bibinfo{year}{2024}) \bibinfo{pages}{133--142}.
\bibitem[{Liu et~al.(2023)Liu, Yao, Peng, and Zhou}]{liu2023bayesian}
\bibinfo{author}{X.~Liu}, \bibinfo{author}{W.~Yao}, \bibinfo{author}{W.~Peng}, \bibinfo{author}{W.~Zhou},
\newblock \bibinfo{title}{Bayesian physics-informed extreme learning machine for forward and inverse pde problems with noisy data},
\newblock \bibinfo{journal}{Neurocomputing} \bibinfo{volume}{549} (\bibinfo{year}{2023}) \bibinfo{pages}{126425}.
\bibitem[{Wang et~al.(2025)Wang, Li, Zhang, Zhou, and Zhou}]{wang2025physics}
\bibinfo{author}{Q.~Wang}, \bibinfo{author}{C.~Li}, \bibinfo{author}{S.~Zhang}, \bibinfo{author}{C.~Zhou}, \bibinfo{author}{Y.~Zhou},
\newblock \bibinfo{title}{Physics-informed extreme learning machine framework for solving linear elasticity mechanics problems},
\newblock \bibinfo{journal}{International Journal of Solids and Structures} \bibinfo{volume}{309} (\bibinfo{year}{2025}) \bibinfo{pages}{113157}.
\bibitem[{Li et~al.(2023)Li, Wu, Huang, Ding, Tai, Liu, and Wang}]{li2023augmented}
\bibinfo{author}{X.~Li}, \bibinfo{author}{J.~Wu}, \bibinfo{author}{Y.~Huang}, \bibinfo{author}{Z.~Ding}, \bibinfo{author}{X.~Tai}, \bibinfo{author}{L.~Liu}, \bibinfo{author}{Y.-G. Wang},
\newblock \bibinfo{title}{Augmented physics informed extreme learning machine to solve the biharmonic equations via fourier expansions},
\newblock \bibinfo{journal}{arXiv preprint arXiv:2310.13947}  (\bibinfo{year}{2023}).
\bibitem[{Yan et~al.(2022)Yan, Vescovini, and Dozio}]{yan2022framework}
\bibinfo{author}{C.~A. Yan}, \bibinfo{author}{R.~Vescovini}, \bibinfo{author}{L.~Dozio},
\newblock \bibinfo{title}{A framework based on physics-informed neural networks and extreme learning for the analysis of composite structures},
\newblock \bibinfo{journal}{Computers \& Structures} \bibinfo{volume}{265} (\bibinfo{year}{2022}) \bibinfo{pages}{106761}.
\bibitem[{Joshi et~al.(2024)Joshi, Snigdha, and Bhattacharya}]{joshi2024physics}
\bibinfo{author}{K.~Joshi}, \bibinfo{author}{V.~Snigdha}, \bibinfo{author}{A.~K. Bhattacharya},
\newblock \bibinfo{title}{Physics informed extreme learning machines with residual variation diminishing scheme for nonlinear problems with discontinuous surfaces},
\newblock \bibinfo{journal}{IEEE Access}  (\bibinfo{year}{2024}).
\bibitem[{Ren et~al.(2025)Ren, Zhuang, Chen, Yu, and Yang}]{ren2025physics}
\bibinfo{author}{F.~Ren}, \bibinfo{author}{P.-Z. Zhuang}, \bibinfo{author}{X.~Chen}, \bibinfo{author}{H.-S. Yu}, \bibinfo{author}{H.~Yang},
\newblock \bibinfo{title}{Physics-informed extreme learning machine (pielm) for stefan problems},
\newblock \bibinfo{journal}{Computer Methods in Applied Mechanics and Engineering} \bibinfo{volume}{441} (\bibinfo{year}{2025}) \bibinfo{pages}{118015}.
\bibitem[{Huang et~al.(2025)Huang, Chen, and Bai}]{huang2025physics}
\bibinfo{author}{L.~Huang}, \bibinfo{author}{L.~Chen}, \bibinfo{author}{R.~Bai},
\newblock \bibinfo{title}{Physics-informed extreme learning machine applied for eigenmode analysis of waveguides and transmission lines},
\newblock \bibinfo{journal}{International Journal of RF and Microwave Computer-Aided Engineering} \bibinfo{volume}{2025} (\bibinfo{year}{2025}) \bibinfo{pages}{6233356}.
\bibitem[{Rao(2007)}]{rao2007vibration}
\bibinfo{author}{S.~S. Rao}, \bibinfo{title}{Vibration of continuous systems}, \bibinfo{publisher}{John Wiley \& Sons}, \bibinfo{address}{Hoboken, NJ}, \bibinfo{year}{2007}.

\end{thebibliography}

\end{document}